\input amstex
\documentstyle{amsppt}
\input xypic
\NoBlackBoxes

\define\DD{\Bbb D}

\define\Q{\Bbb Q}
\define\G{\Bbb G}
\define\Z{\Bbb Z}
\define\M{\Bbb M}
\define\N{\Bbb N}
\define\F{\Bbb F}

\define\R{\Bbb R}
\define\C{\Bbb C}

\define\cal{\Cal}

\document
\topmatter \rightheadtext{On  arithmetic in Mordell-Weil groups}
\title On  arithmetic in Mordell-Weil groups
\endtitle
\author G. Banaszak, P. Kraso{\' n}
\endauthor
\address
Department of Mathematics, Adam Mickiewicz University,
Pozna\'{n}, Poland
\endaddress
\email banaszak\@amu.edu.pl
\endemail
\address
Department of Mathematics, University of Szczecin,
Szczecin, Poland
\endaddress
\email krason\@wmf.univ.szczecin.pl
\endemail

\abstract
In this paper we investigate linear dependence of points
in Mordell-Weil groups of abelian varieties via reduction
maps. In particular we try to determine the conditions
for detecting linear dependence in Mordell-Weil groups via
finite number of reductions.
\endabstract
\endtopmatter

\subhead  1. Introduction
\endsubhead
\medskip

\noindent
The main objective of the paper is to  investigate linear dependence
of points in the Mordell-Weil groups of abelian varieties via the
reduction maps and the height function. In section 5 we prove the following
theorem.

\proclaim{Theorem A} Let $A/F$ be an abelian variety
defined over a number field $F.$ Assume that
$A$ is isogeneous to $A_{1}^{e_1} \times \dots \times A_{t}^{e_t}$
with $A_i$ simple, pairwise nonisogenous abelian varieties
such that $dim_{End_{F^{\prime}}(A_i)^{0}}\,
H_{1} (A_i(\C);\, \Q)\, \geq e_i$
for each $1 \leq i \leq t,$
where $End_{F^{\prime}}(A_{i})^{0} := End_{F^{\prime}}(A_{i})\otimes \Q$
and $F^{\prime}/F$ is a finite extension
such that the isogeny is defined over $F^{\prime}.$
Let $P \in A(F)$ and let
$\Lambda$ be a subgroup of $A(F).$
If $r_v (P) \in r_v (\Lambda)$ for almost all
$v$ of  ${\cal O}_F$ then $P \in \Lambda + A(F)_{tor}.$

\noindent
Moreover if $A(F)_{tor}
\subset \Lambda,$
then the following conditions are equivalent:
\roster
\item[1]\,\, $P \in \Lambda$
\item[2]\,\, $r_v (P) \in r_v (\Lambda)$ for almost all
$v$ of  ${\cal O}_F.$
\endroster
\endproclaim
\noindent

In section 6, Proposition 6.2. we show that the assumption in Theorem A concerning the 
upper bound of the number of simple factors is the best possible in full generality.
The phrase {\it full generality} 
in the previous sentence means {\it for any} $P \in A(F)$ 
{\it and any subgroup} $\Lambda \subset A(F).$
\medskip
 
It has been understood for many years and presented in numerous
papers eg. [Ri] that  many  arithmetic
problems for $\G_m /F$ and methods of treating them  are very similar to those for $A/F.$ This similarity has also been shown in [BGK1] and [BGK2].
Theorem A is an analogue for abelian varieties of a theorem of
A. Schinzel, [Sch, Theorem 2, p. 398],
who proved that for any $\gamma_1,\dots, \gamma_r \in F^{\times}$
and $\beta \in F^{\times}$ such that
$\beta = \prod_{i = 1}^r \gamma_{i}^{n_{v,i}}\,\, \text{mod}\,\, v$
for some $n_{i,v}, \dots, n_{r, v} \in \Z$ and
almost all primes $v$ of ${\cal O}_F$
there are $n_1,\dots, n_r \in \Z$ such that
$\beta = \prod_{i = 1}^r \gamma_{i}^{n_i}.$ The theorem of A. Schinzel
was proved again by Ch. Khare [Kh]
by means of  methods of C. Corralez-Rodrig{\'a}{\~ n}ez and R. Schoof
[C-RS]. The theorem of A. Schinzel concerns the algebraic group
$\G_m / F$ and
does not extend in full generality to
$T = \G_m / F \,\, \times \,\, \G_m /F$ (see [Sch], p. 419). Hence in particular
the theorem of A. Schinzel does not extend in full generality to
algebraic tori and more generally to semiabelian varieties over $F.$
Again the phrase {\it full generality} 
in the last sentence means {\it for any} $P \in T(F)$ 
{\it and any subgroup} $\Lambda \subset T(F).$
In section 6 of this paper we observe that our methods 
of the proof of Theorem A can be used to reprove the A. Schinzel's 
result.
W. Gajda asked a question in 2002 which basically states whether the
analogue of the theorem of Schinzel holds for abelian varieties.
The problem posed by W. Gajda is also called the detecting 
linear dependence problem.

Theorem A strengthens the results of [B], [BGK2], [GG] and [We].
Namely, T. Weston [We] obtained the result stated in Theorem A for
$End_{\overline{F}} (A)$ commutative. 
In [BGK2], together with W. Gajda, we proved
Theorem A for elliptic curves without CM and more generally, for a
class of abelian varieties with $End_{{\overline F}} (A) = \Z ,$
without torsion ambiguity. Moreover we showed, [BGK2] Theorem 2.9,
that for any abelian variety, any free $\cal R$-module $\Lambda
\subset A (F)$ and any $P \in A(F)$ such that $End_{F} (A) \, P$ is a free
$End_{F} (A)$-module the condition (2) of Theorem A implies that there is
$a\in\N $ such that $aP\in\Lambda.$  W. Gajda and K. G{\'
o}rnisiewicz, [GG] Theorem 5.1, showed that the coefficient $a$ in
[BGK2] Theorem 2.9 may be taken to be equal to 1. Very short proof
of Theorem 5.1 of [GG] was also given in [B] Prop. 2.8. 
The main result of [B] states that the problem asked by W. Gajda has 
an affirmative solution for all abelian varieties but with the assumption 
that $End_{F} (A) \, P$ is free $End_{F} (A)$-module and $\Lambda$ is a 
free $\Z$-module which has a
$\Z$-basis linearly independent over $End_{F} (A).$
A. Perucca [Pe2], using methods of [B], [GG] and [Kh], has generalized the
results of [B] and [GG] to the case of a product of an abelian
variety and a torus and removed the assumption in [B] and [GG] that
${\cal R}P$ is a free $\cal R$-module. Recently P. Jossen [Jo] has
given a positive solution to the detecting linear dependence problem
for simple abelian varieties. In his paper P. Jossen uses different methods 
from ours. Due to the A. Schinzel's example [Sch p. 419]
and Proposition 6.2 in this paper, the range of tori and abelian varieties 
for which the detecting linear dependence problem can be solved in full 
generality is 
determined by the example of Schinzel [Sch. p. 419] in the case of tori and 
the upper bound given in our Theorem A in the case of abelian varieties. 
In Section 6, Proposition 6.2 we give an explicit counterexample 
to the problem of detecting linear dependence for the case
of an abelian surface which is a second power of a CM elliptic curve. 
This abelian surface is just beyond our upper bound of Theorem A.  
P. Jossen and A. Perucca [JP] found independently a counterexample 
to the problem of detecting linear dependence.    

\medskip
The proof of Theorem A relies on simultaneous application of
transcendental, $l$-adic and $\pmod  v$ techniques in the theory of
abelian varieties over number fields, use of
semisimplicity of the ring $End (A) \otimes_{\Z} \Q $ and
methods from [B] and [W]. As a corollary of Theorem A one
gets the theorem of T. Weston [W].
\medskip

We  would like to consider a strengthening of Theorem A that could be
used for  computer implementations. With respect to this a natural
problem arises.

\proclaim{Problem} Let $A/F$ be an abelian variety over a number field $F$
and let $P \in A(F)$ and let $\Lambda \subset A(F)$ be a subgroup.
Is there an effectively computable finite set $S^{eff}$ of primes $v$ of
${\cal O}_F,$ depending only on $A,$ $P$ and $\Lambda$
such that the following conditions are equivalent? :
\roster
\item[1]\,\, $P \in \Lambda$
\item[2]\,\, $r_v (P) \in r_v (\Lambda)$ for every
$v \in S^{eff}$
\endroster
\endproclaim

\noindent
We address this problem in section 7. Our main result in this section
is the following theorem.

\proclaim{Theorem B} Let $A/F$ satisfy the hypotheses of Theorem A.
Let $P \in A(F)$ and let $\Lambda$ be a subgroup of $A(F).$
There is a finite set  of primes
$v$ of ${\cal O}_F,$ such that
the condition: $r_v (P) \in r_v (\Lambda)$ for
all $v \in S^{fin}$ implies $P \in \Lambda + A(F)_{tor}.$ \noindent
Moreover if $A(F)_{tor} \subset \Lambda,$ then the following
conditions are equivalent: \roster
\item[1]\,\, $P \in \Lambda$
\item[2]\,\, $r_v (P) \in r_v (\Lambda)$ for
 $v \in S^{fin}.$
\endroster
\endproclaim

In the proof of Theorem B we use the  methods of the proof of
Theorem A, supported by the application of the height pairing
associated with the canonical height function on $A$ [HS], [Sil2] and
the effective Chebotarev's theorem [LO]. The finite set $S^{fin}$
depends on $A, P, \Lambda,$ and the choice of a basis of $c\, A(F),$ 
(see the proof of Theorem 7.7).
\medskip

Important ingredients in the proofs of Theorems 5.1 and 7.7
are Theorems 3.3, 3.6, 7.5 and 7.6 concerning the reduction map.
These theorems refine previous
results of [Bar] and [P] in the case of abelian varieties that are
isogeneous to product of simple, pairwise nonisogeneous abelian varieties.
\bigskip

\subhead  2. Notation and general setup
\endsubhead
\medskip

Let $A/F$ be an abelian variety over a number field $F.$
Let $P, P_1, \dots, P_r \in A(F).$ Put
$\Lambda\, :=\, \sum_{i=1}^r \, \Z\, P_i.$ To prove that
$P \in \sum_{i=1}^r \Z P_i \, +\, T$ in $A (F)$ for some $T \in A(F)_{tor}$
it is enough to prove that $P \in \sum_{i=1}^r \Z P_i \, +\, T^{\prime}$
in $A (L)$ for some finite extension $L / F$ and some
$T^{\prime} \in A(L)_{tor}.$ This is clear since $P, P_1, \dots, P_r \in A(F).$
There is an isogeny $\gamma\, :\, A \rightarrow A_{1}^{e_1}
\times\dots\times A_{t}^{e_t}$
where $A_1, \dots, A_t$ are simple, pairwise nonisogeneous
abelian varieties defined
over certain finite extension $L/F$ and $\gamma$ is also defined over $L.$
To prove that
$P \in \sum_{i=1}^r \Z P_i \, +\, T$ in $A (F)$ for some $T \in A(F)_{tor}$
it is enough to prove that
$\gamma (P) \in \sum_{i=1}^r \Z \gamma (P_i)
+ T^{\prime}$ for some $T^{\prime} \in
\prod_{i=1}^t A_{i}^{e_i} (L)_{tor}.$
Indeed, in this situation there is an element $Q \in \Lambda$ such that
for $M$ equal to the order of $T^{\prime}$ the element
$M (P - Q) \in \text{Ker}\, \gamma.$ Hence
$M (P - Q) \in A (L)_{tor}$ so $(P - Q) \in A (L^{\prime})_{tor}$
where $L^{\prime}/L$ is a finite extension. But $P-Q \in A(F)$ so
$ P \in Q + A(F)_{tor}.$ From now on we can assume that
$A =  A_{1}^{e_1} \times\dots\times A_{t}^{e_t},$ where
$A_1, \dots, A_t$ are simple, pairwise nonisogeneous and defined over $F.$
The remark above shows that we can take $F$ such that
$End_{F} (A_i) = End_{\overline{F}} (A_i)$ for all $i = 1, \dots, t.$
\medskip

We define $r (A) := A_1 \times\dots\times A_t.$
The abelian variety $r (A)$ is called the radical of $A.$
Although it certainly depends on the decomposition of $A$ into simple
factors,  it  is unique up to isogeny.
\medskip

By the remarks above we can assume that
$A =  A_{1}^{e_1} \times\dots\times A_{t}^{e_t}$
where $A_1, \dots, A_t$ are simple abelian varieties defined
over $F.$ Let ${\cal R} := End_{F} (A).$ 
Let ${\cal R}_i := End_{F} (A_i)$ and
$D_i := {\cal R}_i \otimes_{\Z} \Q$ for all $1 \leq i \leq t.$
Then ${\cal R} = \prod_{i=1}^t \, M_{e_i} ({\cal R}_i).$
Let ${\cal L}_i$ be the Riemann lattice such that
$A_i(\C) \cong {\C}^g / {\cal L}_i$ for all $1 \leq i \leq t.$
Then $V_i := {\cal L}_i \otimes_{\Z} \Q$ is a finite dimensional
vector space over $D_i.$ For each $1 \leq i \leq t$ there is a lattice
${\cal L}_{i}^{\prime} \subset {\cal L}_i$ of index
$M_{1,i} := [{\cal L}_i :\, {\cal L}_{i}^{\prime}]$ which is a free
${\cal R}_i$-submodule of ${\cal L}_i$ of rank equal to
$\text{dim}_{D}\,\, V_i .$
Let $K/\Q$ be a finite extension such
that $D_i \otimes_{\Q} K \cong M_{d_i} (K)$ for each $1 \leq i \leq t.$
Hence $V_i \otimes_{\Q} K$ is a free
$M_{d_i} (K)$-module of rank equal to $\text{dim}_{D_i}\,\, V_i .$
Moreover, ${\cal R}_i \otimes_{\Z} {\cal O}_{K} \subset M_{d_i} (K)$
is an ${\cal O}_{K}$ order in $D_i \otimes_{\Q} K \cong M_{d_i} (K)$
and ${\cal L}_{i}^{\prime} \otimes_{\Z} {\cal O}_{K}$ is a free

${\cal R}_i \otimes_{\Z} {\cal O}_{K}$-module
of rank equal to $\text{dim}_{D}\,\, V_i .$
Let $l$ be a prime number. Then $T_l (A_i) \cong
{\cal L}_i \otimes_{\Z} \Z_l$ for every prime number $l \in \Z$
and every $1 \leq i \leq t.$
For a prime ideal  $\lambda | l$ in ${\cal O}_{K}$
let $\epsilon$ denote the index of ramification of $\lambda$
over $l.$
\medskip

Let $L/F$ be a finite extension. From now on $w$ will denote a prime
of ${\cal O}_L$ over a prime $v$ of ${\cal O}_F.$ For a prime
$w$ of good reduction [ST] for $A/L$ let
$$r_w\, :\, A(L) \rightarrow A_w (k_w)$$
be the reduction map.
\medskip

Put $c := |A(F)_{tor}|$ and $\Omega := c\, A(F).$
Note that $\Omega$ is torsion free.
The question we will consider is when the condition
$r_v (P) \in r_v(\Lambda)$ for almost all $v$ of ${\cal O}_{F}$
implies $P \in \Lambda\, + \, A(F)_{tor}.$
\medskip

\noindent
The condition $r_v (P) \in r_v(\Lambda)$ implies 
$r_v (c P) \in r_v(c \Lambda).$ Moreover 
$c\, P \in c\, \Lambda + A(F)_{tor}$
is equivalent to $P \in \Lambda + A(F)_{tor}.$
Hence to answer the question in general it is enough to consider 
the case $P \in c\, A(F),$
$P \not= 0$ and  $\Lambda \subset c \, A(F),$ \, $\Lambda \not= \{0\}.$
\medskip

From now on
we will assume in the proofs of our theorems
that $P \in \Omega,$ \, $P \not= 0,$
$\Lambda \subset \Omega$ and $\Lambda \not= \{0\},$ although
the theorems will be stated for any $P \in A(F)$ and any
subgroup $\Lambda \subset A(F).$
Let $P_1, \dots, P_r,\dots, P_s$ be such a $\Z$-basis of $\Omega$
that:
$$\Lambda = \Z d_1 P_1 + \dots +  \Z d_r P_r + \dots + \Z d_sP_s .
\tag{2.1}$$
where $d_{i} \in \Z \backslash \{0\}$ for $1 \leq i \leq r$ and
$d_i = 0$ for $i > r.$
We put ${\Omega}_{j} := c\, A_j (F).$ Note that
$\Omega = \bigoplus_{j=1}^t\, \Omega_{j}^{e_j}.$
\medskip

For $P \in {\Omega} = {\sum}_{i=1}^{s}{\Bbb Z}P_i$
we write
$$P = n_1 P_1 + \dots + n_r P_r + \dots + n_s P_s \tag{2.2}$$
where $n_i \in \Z.$ Since $\Lambda \subset \Omega$ is a free subgroup of
the free finitely generated abelian group $\Omega,$
observe that $P \in \Lambda$ if and only if  $P \otimes 1 \in
\Lambda \otimes_{\Z} {\cal O}_K.$ The latter is equivalent to
$P \otimes 1 \in \Lambda \otimes_{\Z} {\cal O}_{\lambda}$
for all prime ideals $\lambda \,|\, l$ in ${\cal O}_K$
and all prime numbers $l.$
\bigskip

\noindent
\subhead  3. The reduction map
\endsubhead
\medskip

\noindent
Let $A$ be a product of simple
nonisogenous abelian varieties. Hence
$A = A_1 \times \dots\times A_t$ and in our notation
$e_1 = \dots = e_t = 1.$  To treat this case we need some strengthening
of the results of [BGK2], [Bar] and [P] concerning the reduction
map. Let $L/F$ be any finite extension. Let $P_{i1},\dots ,P_{ir_i} \in
A_i (L)$ be linearly independent over ${\cal R}_i$
for each $1 \leq i \leq t.$ Put $L_{l^{\infty}} := L(A[l^{\infty}]),$
$G_{l^{\infty}} := G({L_{l^{\infty}}/L}),$ \,
$H_{l^{\infty}} := G({\overline F}/ L_{l^{\infty}})$ and
$H_{l^{k}} := G({\overline F}/ L_{l^{k}})$ for all $k \geq 1.$
For each $1 \leq i \leq t$ and $1 \leq j \leq r_i$ let
$$\phi_{ij}\,:\, H_{l^{\infty}} \rightarrow T_{l} (A_i)$$
denote the inverse limit over $k$ of the Kummer maps:
$$\phi_{ij}^{(k)}\,:\, H_{l^{k}} \rightarrow A_i [l^k],$$
$$\phi_{ij}^{(k)} (\sigma) := \sigma\, ({1 \over l^k} P_{ij})
\, -\, {1 \over l^k} P_{ij}.$$

\proclaim{Lemma 3.1 } If ${\alpha}_{11},\dots ,{\alpha}_{1r_1} \in
{\cal R}_1 \otimes_{\Z}{{\Z}}_{l}, \,\,\dots,\,\,
{\alpha}_{t1},\dots ,{\alpha}_{tr_t}
\in {\cal R}_t \otimes_{\Z}{{\Z}}_{l}$ are such that
$\sum_{i=1}^t \sum_{j =1}^{r_t}
{{\alpha}_{ij}}{\phi_{ij}} = 0,$ then \,\,  ${\alpha}_{ij} = 0$ in
${\cal R}_i$ for all $1 \leq i \leq t,$ $1 \leq j \leq r_i.$
\endproclaim

\demo{Proof} Let $\Psi$ be the composition of maps:
$$A(L)\otimes_{\Z}{\Z}_{l}\hookrightarrow
H^{1}(G_{L}; T_{l} (A)) \longrightarrow H^{1}(H_{l^{\infty}};T_{l} (A))
\,\, = \,\,Hom(H_{l^{\infty}};T_{l} (A)).$$
Observe that $\Psi(P_{ij} \otimes 1) = \phi_{ij}.$ By [Se] p. 734
the group $H^1 (G_{l^{\infty}};\, T_{l} (A)$ is finite hence
$ker\Psi\subset (A(L)\otimes_{\Z}{\Z}_{l})_{tor}.$ Let
$c := |A(L)_{tor}| .$ Since $\Psi$ is an ${\cal
R}{\otimes}_{\Z}{{\Z}}_{l}${-}homomorphism, we have:
$$0\quad =\quad
\sum_{i=1}^t \sum_{j =1}^{r_t}
{{\alpha}_{ij}}{\phi_{ij}} \quad =\quad
\Psi(\sum_{i=1}^t \sum_{j =1}^{r_t}
{{\alpha}_{ij}}(P_{ij} \otimes 1)).$$ Hence,
$\sum_{i=1}^t \sum_{j =1}^{r_t}
{{\alpha}_{ij}}(P_{ij} \otimes 1) \in
(A(L)\otimes_{\Z}{\Z}_{l})_{tor}.$ Hence
$ c\, \sum_{i=1}^t \sum_{j =1}^{r_t}
{{\alpha}_{ij}}(P_{ij} \otimes 1)  = 0$ in
$A(L)\otimes_{\Z}{\Z}_{l}.$ Since
$P_{i1}{\otimes}1,\dots ,P_{ir_i} {\otimes}1$ are linearly independent over
${\cal R}_i \otimes_{\Z}{\Z}_{l}$ in $A_i(L)\otimes_{\Z}{\Z}_{l}$ we obtain
$c {{\alpha}_{ij}}=0$ so,
$${{\alpha}_{i1}} = \dots ={{\alpha}_{ir_i}} = 0$$
for each $1 \leq i \leq t$ because ${\cal R}_i$ is
a free ${\Z}${-}module. \qed\enddemo

\noindent
Define the following  maps:
$$\Phi_{i}^k\, :\, H_{l^{k}} \rightarrow A_i [l^k]^{r_i}$$
$$\Phi_{i}^k (\sigma) := (\phi_{i1}^{(k)} (\sigma), \dots,
\phi_{i\, r_i}^{(k)} (\sigma))$$
Then define the map
$\Phi^k\, :\, H_{l^{k}} \rightarrow \bigoplus_{i=1}^{t} A_i[l^k]^{r_i}$
as follows
$\Phi^k := \bigoplus_{i=1}^t \Phi_{i}^k.$
\medskip

\noindent
Define the following maps:
$$\Phi_i\, :\, H_{l^{\infty}} \rightarrow T_l(A_i)^{r_i}$$
$$\Phi_i (\sigma) := (\phi_{i1} (\sigma),
\dots, \phi_{i\, r_i} (\sigma))$$
Again define the map
$\Phi\, :\, H_{l^{\infty}} \rightarrow \bigoplus_{i=1}^{t} T_l(A_i)^{r_i}$
as follows $\Phi := \bigoplus_{i=1}^t \Phi_i.$

\proclaim{Lemma 3.2}The image of the map $\Phi$ is
 open in $\bigoplus_{i=1}^{t}T_l(A_i)^{r_i}.$
\endproclaim

\demo{Proof} Let $T := \bigoplus_{i=1}^{t}T_l(A_i)^{r_i}$ and
let $W := T \otimes_{\Z_l} \Q_l = \bigoplus_{i=1}^{t} V_{i\, l}^{r_i}$
where $ V_{i\, l} := T_l(A_i) \otimes_{\Z_l} \Q_l.$ Denote by ${\Phi}\otimes 1$ the composition of $\Phi$ with the obvious natural inclusion $T\hookrightarrow W.$
Put $M := Im(\Phi\otimes 1) \subset W.$ Both $M$
and $W$ are ${\Q}_{l}[G_{l^{\infty}}]${-}modules.
It is enough to show that $Im\,\Phi$ has a finite index
in $T$ (cf, [Ri, Th. 1.2]). Hence it is enough to show that
$\Phi\otimes 1$ is onto. Observe that $V_{i\, l}$ is a
semisimple ${\Q}_{l}[G_{l^{\infty}}]${-}module for each $1 \leq i \leq t$
because it is a direct summand of the semisimple
${\Q}_{l} [G_{l^{\infty}}]${-}module $V_{l} (A) = \bigoplus_{i=1}^t
V_{i\, l}$ cf. [Fa] Th. 3. Note that $G_{l^{\infty}}$ acts on
$V_{i\, l}$ via the quotient $G (L(A_i[l^{\infty}])/L).$
If $\Phi\otimes 1$ is not onto we have a decomposition
$W = M \oplus M_{1}$ of ${\Q}_{l}[G_{l^{\infty}}]${-}modules
with $M_1$ nontrivial.
Let $\pi_{M_1}\, :\, W \rightarrow W$ be the projection onto
$M_1$ and let $\pi_{i}:W\rightarrow V_{i\, l}$ be a projection that
maps $M_{1}$ nontrivially. Denote
$\widetilde{\pi_i} := \pi_i \, \circ\, \pi_{M_1}.$ By [Fa] Cor 1.
we get $Hom_{G_{l^{\infty}}} (V_{i\, l};\, V_{i^{\prime}\, l})\, \cong\,
Hom_{L} (A_i;\, A_{i^{\prime}}) \otimes_{\Z_l} \Q_l\, =\, 0$
for all $i \not= i^{\prime}.$ Hence
$$\widetilde{\pi_i} (v_{ij})=
\sum_{j=1}^{r_i} {\beta}_{ij} v_{ij},$$
for some ${{\beta}_{ij}}\in {\cal R}_i \otimes \Q_l.$
Since $\pi_{i}$ is nontrivial on $M_1,$
we see that some ${{\beta}_{ij}}$ is
nonzero. On the other hand
$$\widetilde{\pi_i} ({\Phi}(h) \otimes 1) = \sum_{j=1}^{r_i}
{{\beta}_{ij}}({\phi}_{ij}(h)\otimes 1) = 0,$$
for all $h\in H_{l^{\infty}}.$ Since
${{\beta}_{ij}}\in {\cal R}_i \otimes \Q_l,$ we can multiply the
last equality by a suitable power of $l$ to get:
$$0 = \sum_{j=1}^{r_i}
{{\alpha}_{ij}}({\phi}_{ij}(h)\otimes 1),$$
for some
${{\alpha}_{ij}} \in {\cal R}_i\otimes \Z_l.$
Since the maps:
${\cal R}_i \otimes \Z_l \hookrightarrow {\cal R}_i \otimes \Q_l,$
$Hom(H_{l^{\infty}},T_l)\hookrightarrow
Hom(H_{l^{\infty}},V_l)$ are imbeddings of
${\cal R}\otimes\Z_l${-}modules, we obtain \,\,
$\sum_{j=1}^{r_i}{{\alpha}_{ij}}{\phi}_{ij}=0.$
By Lemma 3.1 we get
${{\alpha}_{i1}}=\dots ={{\alpha}_{i\, r_i}}=0,$ hence
${{\beta}_{i1}}=\dots ={{\beta}_{i\, r_i}}=0$ because
${\cal R}$ is torsion free. This contradiction shows
that $M_1=0.$
\qed\enddemo

\proclaim{Theorem 3.3}
Let $Q_{ij} \in A_i (L)$ for $1 \leq j \leq r_i$ be independent over
${\cal R}_i$ for each $1 \leq i \leq t.$
There is a family of primes $w$ of ${\cal O}_L$
of positive density such that
$r_w (Q_{ij}) = 0$ in $A_{i\, w} (k_w)_l$ for all
$1 \leq j \leq r_i$ and $1 \leq i \leq t.$
\endproclaim

\demo{Proof}
\noindent {\bf Step 1.}
We argue in the same way as in the proof of Proposition 2 of [BGK3].
By lemma 3.2 there is an $m \in \N$ such that
$l^m\, \bigoplus_{i=1}^t T_{l} (A_i)^{r_i} \subset
{\Phi}\bigl{(}H_ {l^{\infty}})\bigr{)}
\subset \bigoplus_{i=1}^t T_{l} (A_i)^{r_i}.$
Let $\Gamma$ be the ${\cal R}$-submodule of $A (L)$ generated by all
the points $Q_{ij}.$ Hence $\Gamma =
\sum_{i=1}^t \sum_{j=1}^{r_i} {\cal R}_i Q_{ij}.$ For $k \geq m$ consider
the following commutative diagram.

$$ \CD
G(L_{l^{\infty}}({1 \over l^{\infty}} \Gamma)/L_{l^{\infty}})
@>{\overline{\Phi}}>> \bigoplus_{i=1}^t T_{l} (A_i)^{r_i}/l^{m}
\bigoplus_{i=1}^t T_{l} (A_i)^{r_i}  \\
@VVV @VVV\\
G(L_{l^{k+1}}({1 \over l^{k+1}} \Gamma)/L_{l^{k+1}})
 @>{\overline{\Phi^{k+1}}}>> \bigoplus_{i=1}^t (A_i [l^{k+1}])^{r_i}/
l^{m} \bigoplus_{i=1}^t (A_i [l^{k+1}])^{r_i}\\
@VVV @VV{=}V\\
G(L_{l^{k}}({1 \over l^{k}} \Gamma)/L_{l^{k}})
 @>{\overline{\Phi^k}}>> \bigoplus_{i=1}^t (A_i [l^{k}])^{r_i} /
l^{m} \bigoplus_{i=1}^t (A_i [l^{k}])^{r_i}
\endCD$$

\noindent

The maps $\overline{\Phi}$ and
$\overline{\Phi^{k}},$ for all $k \geq 1,$ are induced
naturally by Kummer maps. For $k \gg 0$ the images
of the middle and bottom horizontal arrows in this diagram are isomorphic.
Hence $G(L_{l^{k+1}}({1 \over l^{k+1}} \Gamma)/L_{l^{k+1}})$ maps onto
$G(L_{l^{k}}({1 \over l^{k}} \Gamma)/L_{l^{k}})$
via the left bottom vertical arrow in the diagram
because the map $\overline{\Phi^k}$ is injective
for each $k \geq 1.$
So quick look at the following tower of fields

$$\diagram
&&&L_{l^{k+1}}(\frac{1}{l^{k+1}}{\Gamma})\dline&\\
&&&L_{l^{k+1}}(\frac{1}{l^{k}}{\Gamma})&\\
&& L_{l^{k}}(\frac{1}{l^{k}}{\Gamma}) \urline&& L_{l^{k+1}} \ulline \\
&&& L_{l^{k}}\ulline_{id}\urline^{h}&&\\
\enddiagram
$$
gives

$$L_{l^k}(\frac{1}{l^k} {\Gamma})\,\,\cap\,\,
L_{l^{k+1}} \quad =\quad L_{l^k} \quad \text{for}\quad k \gg 0
\tag{3.4}$$
\medskip

\noindent {\bf Step 2.} Let $h \in G(L_{l^{\infty}}/L_{l^{k}})$ be
the automorphism which acts on $T_{l}A$ as a homothety $1 + l^k u$ for
some $u \in \Z_{l}^{\times}.$ Such a homothety exists for $k \gg 0$
by the result of Bogomolov [Bo, Cor. 1, p. 702]. Let $h$ also
denote, by a slight abuse of notation, the projection of $h$ onto
$G(L_{l^{k+1}}/L_{l^{k}}).$
By (3.4) we can choose $\sigma \in
G(L_{l^{k+1}}(\frac{1}{l^k} {\Gamma})/ L)$ such that
$\sigma_{|\, L_{l^{k}}(\frac{1}{l^k} {\Gamma})} = \,\, \text{id}$
and $\sigma_{|\,  L_{l^{k+1}}} = h.$
By Chebotarev density theorem there is a
family of primes $w$ of ${\cal O}_L$ of positive density such that
there is a prime $w_1$ in ${\cal O}_{L_{l^{k+1}}(\frac{1}{l^{k}}{\Gamma})}$
over $w$ whose Frobenius in
$L_{l^{k+1}}(\frac{1}{l^{k}}{\Gamma})/L$  equals to $\sigma.$
\medskip

\noindent
Let $l^{c_{ij}}$ be the
order of the element $r_w(Q_{ij})$ in the group
$A_{i\, w} (k_w)_l,$ for some $c_{ij} \geq 0.$
Let $w_2$ be the prime of
${\cal O}_{L_{l^{k}}(\frac{1}{l^k} {\Gamma}))}$
below $w_1.$
Consider the following commutative diagram:
$$ \CD
A_i (L)@> {r_w}>> A_{i\, w} (k_{w})_l\\
@VVV @VV{=}V\\
A_i (L_{l^k}(\frac{1}{l^{k}}{\Gamma}))
@>{r_{w_2}}>> A_{i,\, w} (k_{w_2})_l\\ @VVV @VVV\\
A_i (L_{l^{k{+}1}}(\frac{1}{l^{k}}{\Gamma})
@>{r_{w_1}}>> A_{i\, w} (k_{w_1})_l
 \endCD \tag{3.5}$$

\noindent
Observe that all vertical arrows in the diagram (3.5) are injective.
Let $R_{ij}: = \frac{1}{l^k} Q_{ij} \in
A(L_{l^{k}}(\frac{1}{l^k} \Gamma)) \subset
A(L_{l^{k+1}}(\frac{1}{l^k} \Gamma)).$
The element $r_{w_1}(R_{ij})$ has order $l^{k{+}c_{ij}}$ in
the group $A_{i\, w_1} (k_{w_1})_l$ because
$l^{k{+}c_{ij}} r_{w_{1}} (R_{ij}) = l^{c_{ij}}\, r_w ( Q_{ij}) = 0.$
By the choice of $w,$ we have
$k_w = k_{w_2}$ hence $r_{w_1}(R_{ij})$ comes
from an element of $A_{i\, w}(k_w)_l.$
If $c_{ij} \geq 1$ then
$$h(l^{c_{ij}-1} r_{w_1}( R_{ij})) =
(1 + l^k u) l^{c_{ij}-1} r_{w_1}(R_{ij})$$ since
$l^{c_{ij}-1} r_{w_1}(R_{ij}) \in A_{i\, w}(k_w) [l^{k+1}].$
On the other hand, by
the choice of $w,$ Frobenius at $w_1$ acts on
$l^{c_{ij}-1} r_{w_1}(R_{ij})$ via $h.$ So
$h(l^{c_{ij}-1} r_{w_1}(R_{ij})) = l^{c_{ij}-1} r_{w_1}(R_{ij})$ because
$r_{w_1}(R_{ij}) \in A_{i\, w}(k_w)_l.$
Hence, $l^{c_{ij}-1} u  r_{w_1}(Q_{ij}) = 0$ but this
is impossible since the order of $r_{w_1}(Q_{ij}) = 0$ is $l^{c_{ij}}.$
Hence we must have $c_{ij} = 0.$
\qed\enddemo

\proclaim{Theorem 3.6}
Let $l$ be a prime number.
Let $m \in \N \cup \{0\}$ for all
$1 \leq j \leq r_i$ and $1 \leq i \leq t.$
Let $L/F$ be a finite extension
and let $P_{ij} \in A_i (L)$
be independent over ${\cal R}_i$
and let $T_{ij} \in A_i [l^m]$ be aribitrary torsion elements
for all $1 \leq j \leq r_i$ and $1 \leq i \leq t.$
There is a family of primes $w$ of ${\cal O}_L$
of positive density such that
$$r_{w^{\prime}} (T_{ij}) = r_w (P_{ij}) \,\,\,
\text{in}\,\,\, A_{i, w} (k_w)_l$$
for all $1 \leq j \leq r_i$ and $1 \leq i \leq s,$
where $w^{\prime}$ is a prime in ${\cal O}_{L(A_i [l^m])}$ over $w$
and $r_{w^{\prime}}\, :\, A_{i} (L(A_i [l^m])) \rightarrow
A_{i, w} (k_{w^{\prime}})$ is the corresponding reduction map.
\endproclaim
\demo{Proof} It follows immediately from Theorem 3.3 taking
$L(A[l^m]$ for $L$ and putting $Q_{ij} := P_{ij} - T_{ij}$
for all $1 \leq j \leq r_i$ and $1 \leq i \leq t.$
\qed\enddemo

\proclaim{Remark 3.7}
{\rm  Theorem 3.3 obviously follows from Theorem 3.6.}
\endproclaim
\proclaim{Remark 3.8}
{\rm We have recently learned that A. Perucca  using different methods obtained  analogous theorems to our Theorems 3.3 and 3.6, for the  setting of semiabelian varieties [Pe1 Proposition 11 and 12] . }
\endproclaim

\noindent
\subhead  4. Remarks on semisimple algebras and modules
\endsubhead
\medskip

\noindent
In this section let us recall some basic
properties of modules over semisimple algebras
which will be used in the proof of Theorem 5.1 in the
next section. Let $D$ be a division algebra
and let $K_{i} \subset M_{e} (D)$ denote the
left ideal of $M_{e} (D)$ which consists of $i$-the column
matrices of the form

$$\widetilde{\alpha}_i := \left[ \matrix
0 & \dots a_{1i} & \dots & 0\\
0 & \dots a_{2i} & \dots & 0\\
\vdots  & \vdots & \dots & \vdots\\
0 & \dots a_{e\, i} & \dots & 0
\endmatrix \right] \in K_{i}
$$
Let $W$ be a $D$ vector space and let
$e \in \N$ be a natural number. Then $W^e :=
\underbrace{W \times \dots \times W}_{e-times}$ is a
$M_{e} (D)$-module.
For $\omega \in W$ put
$$\widetilde{\omega} := \left[ \matrix
\omega\\
0\\
\vdots\\
0
\endmatrix \right] \in W^{e},
$$

\proclaim{Lemma 4.1} Every
nonzero simple submodule of the $M_{e} (D)$-module
$W^e$ has the following form
$$K_{1} \widetilde{\omega} = \{\widetilde{\alpha}_1 \widetilde{\omega},\,\,
\widetilde{\alpha}_1 \in K_1\} =
\{\left[ \matrix
a_{11} \, \omega\\
a_{21} \, \omega\\
\vdots\\
a_{e1} \, \omega
\endmatrix \right],\,\, a_{i1} \in D, \,\, 1 \leq i \leq e\}$$
for some $\omega \in W.$
\endproclaim

\demo{Proof} Let $\Delta \subset W^e$ be a simple $M_{e} (D)$-submodule.
Since $M_{e} (D) = \sum_{i=1}^e K_i$ then $\Delta = M_{e} (D)\, \Delta =
\sum_{i=1}^e K_i\, \Delta.$ For each $i,$ \,
$K_i\, \Delta $ is a $M_{e} (D)$-submodule
of $\Delta$  hence $\Delta = K_i\, \Delta$ for some $i$ because
$\Delta$ is simple.
Let $\left[ \matrix
\omega_1\\
\omega_2\\
\vdots\\
\omega_e
\endmatrix \right] \in \Delta$ be a nonzero element.
Again by simplicity of $\Delta$ we obtain
$$\Delta = K_i \Delta \, = \, K_i \left[ \matrix
\omega_1\\
\omega_2\\
\vdots\\
\omega_e
\endmatrix \right] = \{\left[ \matrix
a_{1i}\, \omega_i\\
a_{2i}\, \omega_i\\
\vdots\\
a_{ei}\, \omega_i
\endmatrix \right]: \,\, a_{ji} \in D, \, 1 \leq j \leq e \} \, =
\, K_1 \widetilde{\omega_i}. \qed $$
\enddemo
\medskip

\noindent
Let $e_i \in \N$ and let $D_i$
be a division algebra for each $1 \leq i \leq t.$
We will often use the following notation:
$\DD := \prod_{i=1}^t D_i,$ \,\, ${\text{e}} := (e_1, \dots, e_t)$ \,\,
and \,\, $\M_{\text{e}} (\DD) := \prod_{i=1}^t M_{e_i} (D_i).$
If $W_i$ is a vector space over $D_i$ for each $1 \leq i \leq t$
then the space $W := \bigoplus_{i=1}^t W_{i}^{e_i}$ has a natural
structure of $\M_{\text{e}} (\DD)$-module.

\proclaim{Corollary 4.2}
Every nonzero simple $\M_{\text{e}} (\DD)$-submodule of
$W = \bigoplus_{i=1}^t W_{i}^{e_i}$ has the following form
$$K(j)_{1} \widetilde{\omega (j)} = \{\widetilde{\alpha (j)}_1
\widetilde{\omega (j)}:\,\,\,
\widetilde{\alpha (j)}_1 \in K(j)_1\} =
\{\left[ \matrix
a_{11} \, \omega (j)\\
a_{21} \, \omega (j)\\
\vdots\\
a_{e_{j}1} \, \omega (j)
\endmatrix \right],\,\, a_{k1} \in D_j, \,\, 1 \leq k \leq e_{j}\}$$
for some $1 \leq j \leq t$ and some $\omega (j) \in W_j$ where
$K(j)_{1} \subset M_{e_j} (D_j)$ denotes the
left ideal of $M_{e_j} (D_j)$ which consists of $1st$ column
matrices.
\endproclaim

\demo{Proof} Follows immediately from Lemma 4.1. \qed
\enddemo
\bigskip

\noindent
Let $D_i$ be a finite dimensional division algebra over $\Q$
for every $1 \leq i \leq t.$
The trace homomorphisms:
$tr_i \, :\, M_{e_i} (D_i) \rightarrow \Q,$
for all $1 \leq i \leq t,$ give the trace homomorphism
$tr \, :\, \M_{\text{e}} (\DD) \rightarrow \Q,$ where
$tr := \sum_{i=1}^t \, tr_i.$
Let $W_i$ be a finite dimensional $D_i$-vector space
for each $1 \leq i \leq t.$
Then $W$ is a finitely
generated $\M_{\text{e}} (\DD)$-module.
The homomorphism $tr$ gives a natural map of $\Q$-vector spaces
$$ tr \, : \, Hom_{\M_{\text{e}} (\DD)}
(W, \,\, \M_{\text{e}} (\DD))
\rightarrow Hom_{\Q} (W ,\,\, \Q).
\tag{4.3}$$

\proclaim{Lemma 4.4} The map (4.3) is an isomorphism.
\endproclaim
\demo{Proof}
For each $1 \leq i \leq t$ we have the trace map
$$ tr_i \, : \, Hom_{M_{e_i} (D_i)}
(W_{i}^{e_i}, \,\, M_{e_i} (D_i))
\rightarrow Hom_{\Q} (W_{i}^{e_i} ,\,\, \Q).
\tag{4.4}$$
The map (4.3) is naturally compatible with maps (4.4) via
natural isomorphisms:
$$\bigoplus_{i=1}^t \,\, Hom_{M_{e_i} (D_i)}
(W_{i}^{e_i}, \,\, M_{e_i} (D_i)) \cong
Hom_{\M_{\text{e}} (\DD)}  (W, \,\, \M_{\text{e}} (\DD))
\tag{4.5}$$
$$\bigoplus_{i=1}^t \,\,
Hom_{\Q}
(W_{i}^{e_i}, \,\, \Q) \cong Hom_{\Q} (W ,\,\, \Q)
\tag{4.6}$$
In other words $tr = \sum_{i=1}^t \, tr_i.$ Hence it is enough to
prove the lemma for the maps (4.4).
Since $M_{e_{i}}(D_{i})$ is a simple ring for which every simple module
is isomorphic to $K(i)_{1}$ it is enough
to prove that
$$tr_{i}: Hom_{M_{e_{i}}(D_{i})}(K(i)_{1}; M_{e_{i}}(D_{i})) \,\, {\cong}
\,\, Hom_{\Bbb Q}(K(i)_{1}; {\Bbb Q}). \tag{4.7} $$

\noindent
Notice that every map
${\phi}\in Hom_{M_{e_{i}}(D_{i})}(K(i)_{1}; M_{e_{i}}(D_{i}))$ is
determined by its image on the element
$\left [\matrix
1 & 0 & \dots & 0\\
0 & 0 & \dots & 0 \\
\vdots & \vdots & \vdots & \vdots \\
0 & 0 & \dots & 0 \\
\endmatrix \right ] \in K(i)_{1}.$
Since $\phi$ is a $M_{e_{i}}(D_{i})$-module homomorphism we have

$${\phi} \bigl( \left [\matrix
1 & 0 & \dots & 0\\
0 & 0 & \dots & 0 \\
\vdots & \vdots & \vdots & \vdots \\
0  & 0 & \dots &0 \\
\endmatrix \right ] \bigr) =
\left [\matrix
c_{11} & c_{12} & \dots & c_{1\, e_{i}}\\
0 & 0 & \dots & 0 \\
\vdots & \vdots & \vdots & \vdots \\
0 & 0 & \vdots & 0 \\
\endmatrix \right ] \tag{4.8}$$
for some $c_{11}, c_{12}, \dots, c_{1\, e_i} \in D_i.$
The map (4.7) is injective (cf. [Re] Theorem 9.9 ).
From the definition of $K(i)_{1}$ and (4.8) it follows that
dimensions of the $\Q$-vector spaces
$Hom_{M_{e_{i}}(D_{i})}(K(i)_{1}; M_{e_{i}}(D_{i}))$ and
$Hom_{\Q}(K(i)_{1}; {\Q})$ are equal. Hence (4.7) is an isomorphism.
\qed
\enddemo
\medskip

\noindent
The algebra $\M_{\text{e}} (\DD)$ is semisimple hence
the module $W$ is semisimple so for every
$\tilde{\pi} \in
Hom_{\M_{\text{e}} (\DD)}  (W, \,\, \M_{\text{e}} (\DD)) $ there
is a $\M_{\text{e}}(\DD)$-homomorphism
$\tilde{s} : \text{Im}\, \widetilde{\pi} \rightarrow W $ such that
$\widetilde{\pi} \circ \widetilde{s} = Id.$ Because of (4.5) we can write
$\widetilde{\pi} = \prod_{i=1}^t \widetilde{\pi (i)}$ for some
$\widetilde{\pi (i)} \in
Hom_{M_{e_i} (D_i)}  (W_{i}^{e_i}, \,\, M_{e_i} (D_i)).$
Note that $ \text{Im}\, \widetilde{\pi} \, = \, \prod_{i=1}^t\,
\text{Im} \, \widetilde{\pi (i)}.$ For each $1 \leq i \leq t$ we can find
$M_{e_i}(D_i)$-homomorphism
$\widetilde{s (i)} : \,\text{Im} \,  \widetilde{\pi (i)} 
\rightarrow W_{i}^{e_i} $
such that $\widetilde{\pi (i)} \circ \widetilde{s (i)} = Id$ and
$\widetilde{s} = \bigoplus_{i=1}^t \, \widetilde{s (i)}$
because $M_{e_i} (D_i)$ is simple.
\medskip

\noindent
By [Re], Theorem 7.3 every simple $M_{e_i} (D_i)$-submodule of
$M_{e_i} (D_i)$ is isomorphic to $K(i)_{1}.$
Since  $dim_{D_i}\, M_{e_i} (D_i) = e_{i}^2$ and
$dim_{D_i}\, K(i)_{1} = e_{i}$ we see that
$M_{e_i} (D_i)$ is a direct sum of $e_i$ simple $M_{e_i} (D_i)$-submodules.
Hence every $M_{e_i} (D_i)$-submodule of
$M_{e_i} (D_i)$ is a direct sum
of at most $e_i$ simple $M_{e_i} (D_i)$-submodules.

\subhead  5. Detecting linear dependence in Mordell-Weil groups
\endsubhead
\medskip

\noindent
\proclaim{Theorem 5.1}
Let $A/F$ be an abelian variety
defined over a number field $F.$ Assume that
$A$ is isogeneous to $A_{1}^{e_1} \times \dots \times A_{t}^{e_t}$
with $A_i$ simple, pairwise nonisogenous abelian varieties
such that $dim_{End_{F^{\prime}}(A_i)^{0}}\,
H_{1} (A_i(\C);\, \Q)\, \geq e_i$
for each $1 \leq i \leq t$
and $F^{\prime}/F$ is a finite extension
such that the isogeny is defined over $F^{\prime}.$
Let $P \in A(F)$ and let
$\Lambda$ be a subgroup of $A(F).$
If $r_v (P) \in r_v (\Lambda)$ for almost all
$v$ of  ${\cal O}_F$ then $P \in \Lambda + A(F)_{tor}.$

\noindent
Moreover if $A(F)_{tor}
\subset \Lambda,$
then the following conditions are equivalent:
\roster
\item[1]\,\, $P \in \Lambda$
\item[2]\,\, $r_v (P) \in r_v (\Lambda)$ for almost all
$v$ of  ${\cal O}_F.$
\endroster
\endproclaim
\demo{Proof} Assume that $P \notin \Lambda.$ This implies that
$P \otimes 1 \notin \Lambda \otimes_{\Z}
{\cal O}_{\lambda}$ for some $\lambda \,|\, l$ for some prime number
$l.$ Hence in (2.2) $n_j \not= 0$ for some $1 \leq
j \leq s.$ We can consider the equality (2.2) in $\Omega
\otimes_{\Z} {\cal O}_K.$
Since $P \notin \Lambda \otimes_{\Z} {\cal O}_{\lambda}$ then
there is $1 \leq j_0 \leq s$ such that $\lambda^{m_1}\, ||\, n_{j_0}$ and
$\lambda^{m_2}\, |\, d_{j_0}$ for natural numbers $m_1 < m_2.$
Consider the map of $\Z$-modules
$$\pi\, :\, \Omega \rightarrow \Z$$
$$\pi (R) :=  {\mu}_{j_{0}}$$
for $R = \sum_{i=1}^s {\mu}_i P_i $ with $\mu_{i} \in \Z$
for all $1 \leq i \leq s.$ By abuse of notation denote also by $\pi$
the map $\pi \otimes \Q \, :\, \Omega \otimes_{\Z} \Q\rightarrow \Q.$
By Lemma 4.4 there is map $\widetilde{\pi} \in  Hom_{\M_{\text{e}} (\DD)}
(\Omega \otimes_{\Z} \Q, \,\, \M_{\text{e}} (\DD))$
such that $tr( \widetilde{\pi}) = \pi.$ By remarks after proof of Lemma 4.4
there is $\widetilde{s} \in Hom_{\M_{\text{e}} (\DD)}
(\text{Im} \widetilde{\pi}  , \,\, \Omega \otimes_{\Z} \Q)$
such that $\widetilde{\pi} \circ \widetilde{s} = Id.$ Moreover for all
$1 \leq i \leq t$ there are
$\widetilde{\pi (i)} \in  Hom_{M_{e_i} (D_i)}
(\Omega_{i}^{e_i} \otimes_{\Z} \Q, \,\, M_{e_i} (D_i))$ and
$\widetilde{s (i)} \in Hom_{M_{e_i} (D_i)}
(\text{Im} \widetilde{\pi (i)}  , \,\, \Omega_{i}^{e_i} \otimes_{\Z} \Q)$
such that $\widetilde{\pi (i)} \circ \widetilde{s (i)} = Id$ and
$\widetilde{\pi} = \prod_{i=1}^t \widetilde{\pi (i)},$ \,
$\widetilde{s} = \prod_{i=1}^t \widetilde{s (i)}.$ Moreover
$ \, \text{Ker} \, \widetilde{\pi} =
\prod_{i}^t \,\, \text{Ker} \,\, \widetilde{\pi (i)}$ and we have
$\Omega_{i}^{e_i} \otimes_{\Z} \Q \cong  \, \text{Im} \, \widetilde{s (i)}
 \, \oplus \, \text{Ker} \, \widetilde{\pi (i)}$ and
$\Omega \otimes_{\Z} \Q \cong \text{Im} \, \widetilde{s} \oplus
\text{Ker} \, \widetilde{\pi}.$ By Lemma 4.1 we can present
$\text{Im} \, \widetilde{s (i)}$ and
$\text{Ker} \, \widetilde{\pi (i)}$ as direct sums of
simple $M_{e_i} (D_i)$-submodules as follows:
$$\text{Im} \, \widetilde{s (i)} = \bigoplus_{k =1}^{k_i}\,
K(i)_{1} \, \widetilde{\omega_{k} (i)},$$
$$\text{Ker} \, \widetilde{\pi (i)} = \bigoplus_{k = k_{i} + 1}^{u_i}\,
K(i)_{1} \, \widetilde{\omega_{k}(i)}.$$
Observe that $k_i \leq e_i$ for every $1 \leq i \leq t.$
It is simple to observe that
the elements $\omega_{1}(i), \dots, \omega_{k_i}(i), \dots,
\omega_{u_i}(i)$ give a basis of the $D_i$-vector space
$\Omega_i \otimes_{\Z} \Q.$ We can assume without loss of generality that
$\omega_{k_{i} + 1}(i), \dots,
\omega_{u_i}(i) \in \Omega_i.$ Tensoring the map $\pi$ with ${\cal O}_K$
we will denote the resulting map
$\pi\, :\, \Omega \otimes_{\Z} {\cal O}_{K} \rightarrow {\cal O}_{K}$
also by $\pi.$ Similarly tensoring the maps
$\widetilde{\pi (i)}$ and $\widetilde{s (i)}$ with $K$ we get
$M_{e_i} (D_i) \otimes_{\Q} K$-linear homomorphisms
$\widetilde{\pi (i)} \, :\, \Omega_{i}^{e_i} \otimes_{\Z} K \,\,
\rightarrow \,\, M_{e_i} (D_i) \otimes_{\Q} K$ and
$\widetilde{s (i)} \, :\,\, \text{Im} \tilde{\pi}_i
\rightarrow \Omega_{i}^{e_i} \otimes_{\Z} K$ also denoted by
$\widetilde{\pi (i)}$ and $\widetilde{s (i)}$ respectively. Note that for each
$1 \leq i \leq t$ the $K$-vector space $\Omega_i \otimes_{\Z} K$
is a free $D_i \otimes_{\Q} K \cong M_{d_i} (K)$ module.
Recall that ${\cal R} \subset
\M_{\text{e}} (\DD),$ ${\cal R} \otimes_{\Z} \Q =
\M_{\text{e}} (\DD)$ and $\Omega$ is a finitely
generated ${\cal R}$-module.
Hence there is a natural number $M_0$ such that the
homomorphisms of ${\cal R} \otimes_{\Z} {\cal O}_K$-modules are well
defined:
$$M_0\, \widetilde{\pi}\,:\, \Omega \otimes_{\Z} {\cal O}_K
\rightarrow {\cal R} \otimes_{\Z} {\cal O}_K,$$
$$\widetilde{s}\, :\, M_{0}\, \widetilde{\pi} (\Omega \otimes_{\Z} {\cal O}_K)
\rightarrow \Omega \otimes_{\Z} {\cal O}_K,$$
We can restrict the trace homomorphism
to ${\cal R} \otimes_{\Z} {\cal O}_K \subset D \otimes_{\Q} K$
to get an ${\cal O}_K$-linear homomorphism
$tr \, :\, {\cal R} \otimes_{\Z} {\cal O}_K \rightarrow K.$
Note that $tr \, M_0 \, \tilde{\pi} = M_0 \, \pi$ and
$M_0 \, \widetilde{\pi} \circ \widetilde{s} =
M_0 \, Id_{M_{0} \, \tilde{\pi} (\Omega \otimes_{\Z} {\cal O}_K)}.$
Consider now the first column vectors
$K(i)_1 \subset M_{e_i}({\cal R}_i \otimes_{\Z} {\cal O}_K).$
Define the $M_{e_i}({\cal R}_i \otimes_{\Z} {\cal O}_K)$-module
$$\widetilde{{\Gamma} (i)} := \sum_{k=1}^{k_i}\,\,
K(i)_1\,\, M_0 \,\,   \widetilde{\omega_k (i)}
+ \sum_{k=k_{i}+1}^{u_i} K(i)_1 \,\,
 \widetilde{\omega_k (i)} \subset {\Omega}_{i}^{e_i} \otimes_{\Z} {\cal O}_K$$
and ${\cal R} \otimes_{\Z} {\cal O}_K$-module
$\widetilde{{\Gamma}} :=
\bigoplus \widetilde{{\Gamma} (i)} \subset {\Omega} \otimes_{\Z} {\cal O}_K.$
Put $M_2 := [\Omega \otimes_{\Z} {\cal O}_{K}\, :\, \widetilde{\Gamma}]$
and  $M_3 :=
[\widetilde{\Gamma}\, :\, M_2\, \Omega \otimes_{\Z} {\cal O}_{K}].$
By the choice of the point $P_{j_0}$ we get
$\pi (P) \notin \pi (\Lambda \otimes_{\Z} {\cal O}_{\lambda}) +
\lambda^m \,
\pi (\Omega \otimes_{\Z} {\cal O}_{\lambda})$ for $m > m_2.$
Hence
$$M_0\, \widetilde{\pi} (P) \notin
M_0\, \widetilde{\pi} (\Lambda \otimes_{\Z} {\cal O}_{\lambda}) +
M_0\, \lambda^m \, \widetilde{\pi} (\Omega \otimes_{\Z} {\cal O}_{\lambda})
\tag{5.2}$$
because $tr M_0 \widetilde{\pi} = M_0 \pi.$
Put $K(i)_{1, \lambda}  := K(i)_{1} \otimes_{{\cal O}_K} {\cal O}_{\lambda}
\subset M_{e_i}({\cal R}_{i\, \lambda}).$
Let $Q \in \Lambda$ be an arbitrary element.
We can write
$$M_2 P = \sum_{i=1}^t  \sum_{k=1}^{k_i}\,\,
\widetilde{\alpha_k (i)}_{1} M_0\,
\widetilde{\omega_k (i)} +
\sum_{i=1}^t\sum_{k=k_{i}+1}^{u_i}  \widetilde{\alpha_k (i)}_{1}
 \widetilde{\omega_k (i)},$$
$$M_2 Q = \sum_{i=1}^t  \sum_{k=1}^{k_i} \,\, \widetilde{\beta_k (i)}_{1} M_0\,
\widetilde{\omega_k (i)} +
\sum_{i=1}^t\sum_{k=k_{i}+1}^{u_i}  \widetilde{\beta_k (i)}_{1}
 \widetilde{\omega_k (i)},$$
for some $ \widetilde{\alpha_k (i)}_{1},\,
 \widetilde{\beta_k (i)}_{1} \in {K(i)_{1, \lambda}}$
with $1 \leq k \leq u_i$ and $1 \leq i \leq t.$
Then
$$M_0\, \tilde{\pi} (M_2 (P - Q)) =
M_{0}^2 \, \prod_{i=1}^t  \sum_{k=1}^{k_i}\,\,
(\widetilde{\alpha_k (i)}_{1} - \widetilde{\beta_k (i)}_{1}) \,
\widetilde{\pi} (\widetilde{\omega_k (i)}).  \tag{5.3}$$
Since $\widetilde{\pi} = \prod_{i=1}^t \widetilde{\pi (i)}$ maps the 
module $\Omega \otimes_{\Z} \Q = \bigoplus_{i=1}^t 
\Omega_{i}^{e_i} \otimes_{\Z} \Q $ into the ring 
$\M_{\text{e}} (\DD) = \prod_{i=1}^t M_{e_i} (D_i)$ componentwise, 
we replaced $\sum_{i=1}^t$ by $\prod_{i=1}^t.$
Hence (5.2) and (5.3) give
$M_{0}^2 \, \prod_{i=1}^t  \sum_{k=1}^{k_i}\,\,
(\widetilde{\alpha_k (i)}_{1} - \widetilde{\beta_k (i)}_{1}) \,
\widetilde{\pi} (\widetilde{\omega_k (i)})
 \, \notin
\lambda^m M_0\, \tilde{\pi} (M_2\,
\Omega \otimes_{\Z} {\cal O}_{\lambda}),$
so
$$M_{0}^2 \, \prod_{i=1}^t  \sum_{k=1}^{k_i}\,\,
(\widetilde{\alpha_k (i)}_{1} - \widetilde{\beta_k (i)}_{1}) \,
\widetilde{\pi} (\widetilde{\omega_k (i)})
 \,  \notin
\lambda^m M_0 \, \tilde{\pi} (M_3\, \widetilde{\Gamma}). \tag{5.4}$$
Hence for some $1 \leq i \leq t$ and $1 \leq k \leq k_i$ we obtain
$$\widetilde{\alpha_k (i)}_{1} - \widetilde{\beta_k (i)_{1}} \,\notin
\lambda^m M_3 \, K(i)_{1, \lambda}\, . \tag{5.5}$$
Let $\epsilon \in \N$ be the ramification index of
$\lambda$ over $l.$
Observe that for every $n \in \N$ we have an isomorphism
$A_i [\lambda^{\epsilon n}] \cong {\cal L_i} \otimes_{\Z}
{\cal O}_{\lambda}\, /\, \lambda^{\epsilon n}\, {\cal L}_i \otimes_{\Z}
{\cal O}_{\lambda}$ because $l\, {\cal O}_K =
\prod_{\lambda\, |\, l} \lambda^{\epsilon},$
$A_i [l^n] \cong {\cal L}_i \otimes_{\Z} {\Z}_{l}\, /\,
l^n\, {\cal L}_i \otimes_{\Z} {\Z}_{l}$
and $A_i [l^n] = \bigoplus_{\lambda\, |\, l} A_i [\lambda^{\epsilon \, n}].$
Recall that we chose, for each $1 \leq i \leq t,$ a lattice
${\cal L}_{i}^{\prime} \subset {\cal L}_{i}$ such that
${\cal L}_{i}^{\prime}$ is a free ${\cal R}_{i}$-module.
Let $M_4 := \text{max}_{1 \leq i \leq t}
[{\cal L}_i:\, {\cal L}_{i}^{\prime}].$
Put ${\cal L} := \bigoplus_{i=1}^t\, {\cal L}_{i}$ and
${\cal L}^{\prime} := \bigoplus_{i=1}^t {\cal L}_{i}^{\prime}.$
By Snake Lemma
the kernel of the following natural map of ${\cal O}_{\lambda}$-modules
is finite and annihilated by $\lambda^{\epsilon \, m_4}$
$$z(n, \lambda)\, :\, {\cal L}^{\prime} \otimes_{\Z} {\cal O}_{\lambda}\, /\,
\lambda^{\epsilon \, n}\, {\cal L}^{\prime} \otimes_{\Z}
{\cal O}_{\lambda} \rightarrow {\cal L} \otimes_{\Z} {\cal O}_{\lambda}\, /\,
\lambda^{\epsilon \, n}\, {\cal L}
\otimes_{\Z} {\cal O}_{\lambda}, \tag{5.6} $$
where $l^{m_4}\, ||\, M_4.$ Let $m_0$ and $m_3$ denote
the natural numbers with the property
$l^{m_0}\, ||\, M_0$ and $l^{m_3}\, ||\, M_3.$
Let ${\eta}_{1}(i), \dots, {\eta}_{p_i} (i)$ be a basis of
${\cal L}_{i}^{\prime}$ over ${\cal R}_{i}.$ By the assumptions
$p_i \geq e_i.$ Hence ${\cal L}_{i}^{\prime} \otimes_{\Z}
{\cal O}_{\lambda}\, /\, \lambda^{\epsilon \, n}\,
{\cal L}_{i}^{\prime} \otimes_{\Z} {\cal O}_{\lambda}$ is a free
${\cal R}_{i, \lambda}\,/\, \lambda^{\epsilon\, n}\,
{\cal R}_{i, \lambda}$-module with basis
$\overline{{\eta}_{1}(i)}, \dots, \overline{{\eta}_{p_i}(i)},$ where
$\overline{{\eta}_{k}(i)}$ denotes the image of ${\eta}_{k}(i)$ in
${\cal L}^{\prime} \otimes_{\Z} {\cal O}_{\lambda}\, /\,
\lambda^{\epsilon \, n}\, {\cal L}^{\prime} \otimes_{\Z}
{\cal O}_{\lambda}$ for
each $1 \leq k \leq p_i.$ Let $T_{k} (i)$ be the image of
$\overline{{\eta}_{k}(i)}$ via the map
$z(n, \lambda)$ for all $1 \leq i \leq t$ and $1 \leq k \leq p_i.$
Take $n \in \N$ such that $\epsilon \, n > m + \epsilon\, m_0 + \epsilon \,
m_3 + \epsilon \, m_4$
and put $L := F(A[l^{n}]) = F(r(A) [l^{n}]).$
Observe that $A[l^{n}] \subset A(L).$ By Theorem 3.6 there is a
family of primes $w$ of ${\cal O}_L$ of positive density such that
$r_w (\omega_k (i)) = 0$ for $1 \leq i \leq t,$ $k_{i}+1 \leq k \leq u_i$
and $r_w (\omega_k (i)) = r_w (T_{k} (i))$ for all $1 \leq i \leq t,$
$1 \leq k \leq k_i.$
Since $r_w (P) \in r_w (\Lambda)$ we take $Q \in \Lambda$ such that
$r_w (P) = r_w (Q).$ Applying the reduction map $r_w$ to the equation

$$M_2 (P - Q) = \sum_{i=1}^t  \sum_{k=1}^{k_i}\,\,
(\widetilde{\alpha_k (i)}_{1} - \widetilde{\beta_k (i)}_{1}) M_0\,
\widetilde{\omega_k (i)} +
\sum_{i=1}^t\sum_{k=k_{i}+1}^{u_i}
(\widetilde{\alpha_k (i)}_{1} -\widetilde{\beta_k (i)}_{1})
 \widetilde{\omega_k (i)},$$
we obtain
$$0 =
\sum_{i=1}^t  \sum_{k=1}^{k_i}\,\,
(\widetilde{\alpha_k (i)}_{1} - \widetilde{\beta_k (i)}_{1}) M_0\,
\widetilde{r_w (T_k (i))}.$$
Since the map $r_w$ is injective on
$l$-torsion subgroup of $A(L)$ ([HS] Theorem C.1.4 p. 263,
[K] p. 501-502),
we obtain
$$0 =
\sum_{i=1}^t  \sum_{k=1}^{k_i}\,\,
(\widetilde{\alpha_k (i)}_{1} - \widetilde{\beta_k (i)}_{1}) M_0\,
\widetilde{T_k (i)}.$$
Therefore $\sum_{i=1}^t  \sum_{k=1}^{k_i}\,\,
(\widetilde{\alpha_k (i)}_{1} - \widetilde{\beta_k (i)}_{1}) M_0\,
\widetilde{\overline{\eta_{k} (i)}} \, \, \in \,\,
\text{Ker}\,\, z(n, \lambda).$
So, the element
\noindent
 $\lambda^{\epsilon\, m_0 + \epsilon\, m_4}
\sum_{i=1}^t  \sum_{k=1}^{k_i}\,\,
(\widetilde{\alpha_k (i)}_{1} - \widetilde{\beta_k (i)}_{1}) \,
\widetilde{\eta_{k} (i)}$ maps to zero in
${\cal L}^{\prime} \otimes_{\Z} {\cal O}_{\lambda}\, /\,
\lambda^{\epsilon \, n}\, {\cal L}^{\prime} \otimes_{\Z}
{\cal O}_{\lambda}.$ Hence
$$\sum_{i=1}^t  \sum_{k=1}^{k_i}\,\,
(\widetilde{\alpha_k (i)}_{1} - \widetilde{\beta_k (i)}_{1}) \,
\widetilde{\eta_{k} (i)} \in \lambda^{\epsilon n - \epsilon \, m_0
- \epsilon \, m_4}
\, {\cal L}^{\prime} \otimes_{\Z} {\cal O}_{\lambda}.$$
Since $\eta_{1} (i), \dots, \eta_{p_i} (i)$ is a basis of
${\cal L}_{i}^{\prime} \otimes_{\Z} {\cal O}_{\lambda}$ over
${\cal R}_{i, \lambda},$ we obtain
$$ \widetilde{\alpha_k (i)}_{1} - \widetilde{\beta_k (i)}_{1} \, \in \,
\lambda^{\epsilon n - \epsilon \, m_0 - \epsilon
\, m_4}\, K (i)_{1, \lambda} \tag{5.7}$$
for all $1 \leq i \leq t$ and $1 \leq k \leq k_i.$
But (5.7) contradicts  (5.5) because we chose $n$ such that
$\epsilon\, n - \epsilon\, m_0 - \epsilon \, m_4 > m + \epsilon \, m_3.$
\qed\enddemo

\proclaim{Corollary 5.8} (Weston [We p. 77])
Let $A$ be an abelian variety
defined over a number field such that $End_{\overline{F}} (A)$
is commutative. Then Theorem 5.1 holds for $A.$
\endproclaim

\demo{Proof}
Since $End_{\overline{F}} (A)$ is commutative,
$A$ is isogeneous to $A_{1} \times \dots \times A_{t}$
with $A_i$ simple, pairwise nonisogenous.
In this case the assumption in Theorem 5.1
$dim_{End_{F^{\prime}}(A_i)^{0}}\,
H_{1} (A_i(\C);\, \Q)\, \geq 1$
for each $1 \leq i \leq t$ always holds.
\qed\enddemo

\proclaim{Corollary 5.9} Let $A = E_{1}^{e_1} \times \dots
\times E_{t}^{e_t},$ where $E_{1}, \dots, E_{t}$ are pairwise
nonisogenous elliptic curves defined over $F.$ Assume that  
$1 \leq e_i \leq 2$ if
$End_{F} (E_i) = \Z$  and $e_{i}=1$ if
$End_{F} (E_i) \neq\Z .$ Then Theorem 5.1 holds for $A.$
\endproclaim

\demo{Proof} Observe that for an elliptic curve $E/F$ we have
 $dim_{End_{F}(E)^{0}}\,
H_{1} (E(\C);\, \Q)\, =2$ if $End_{F} (E)= \Z$ and
$dim_{End_{F} (E)^{0}}\, H_{1} (E(\C);\, \Q)\, =1$ if
$End_{F} (E)\neq \Z$
 \qed\enddemo

\proclaim{Remark 5.10}
{\rm Theorem 5.1 and in particular Corollary 5.9 answer the question of
T. Weston [We] p. 77 concerning the noncommutative
endomorphism algebra case.}
\endproclaim

\bigskip
\subhead 6. Counterexamples to the problem of detecting linear dependence via reduction maps
\endsubhead
\bigskip

The hypothesis in Theorem 5.1 that $A$ is isogeneous over $F^{\prime}$ to $A_{1}^{e_1} \times \dots \times A_{t}^{e_t}$ with $dim_{End_{F^{\prime}}(A_i)^{0}}\,
H_{1} (A_i(\C);\, \Q)\, \geq e_i$ for each $1 \leq i \leq t,$ cannot be omitted
in full generality. In fact in Proposition 6.2 we produce counterexamples for the problem of detecting linear dependence, when the hypotheses in Theorem 5.1 does not hold, 
considering products of two CM elliptic curves. 
In this way we show  that the upper bound condition for the number of simple factors in Theorem 5.1 is the best possible as far as full generality is concerned. The idea of the proof of Proposition 6.2 that our family of  abelian  varieties provides counterexamples to Theorem 5.1 is based on the counterexample of A. Schinzel [Sch p.419] for the product of two $\G_m$. 
For this reason let us start the discussion of counterexamples for algebraic tori applying
[Sch p.419]. 

\medskip
\noindent\bf\it The case of algebraic tori

\rm\medskip

\noindent
Let us mention that the methods of the proof of Theorem 5.1 work
for some algebraic tori over a number field $F.$ 
To understand for which 
tori our methods work let $T/F$ be an algebraic torus and 
let $F^{\prime}/F$ be a finite extension
that splits $T.$ Hence $T \otimes_{F} F^{\prime} \cong 
\G_{m}^e := \underbrace{\G_{m} \times \dots\times \G_{m}}_{e-times}$
where $\G_{m} := \,\, \text{spec}\,\, F^{\prime}[t,\, t^{-1}].$
For any field extension $F^{\prime} \subset M \subset \overline{F}$ we have 
$End_{M}\,\, (\G_m) = \Z$ and $H_1 (\G_m (\C);\, \Z) = \Z.$
Hence the condition $e \leq dim_{End_{F^{\prime}}\,\, (\G_m)^{0}}\,\, 
H_1 (\G_m (\C);\, \Q) = 1,$ analogous to the corresponding condition 
of Theorem 5.1, means that $e =1.$ Hence we can prove
the analogue of Theorem 5.1 for one dimensional tori which is basically
the A. Schinzel's Theorem 2 of [Sch]. 
Observe that torsion ambiguity that appears in Theorem 5.1 can be 
removed in the case of one dimensional tori 
by use of an argument similar to the proof of 
Theorem 3.12 of [BGK2]. 
On the other hand A. Schinzel showed that his theorem does not extend in 
full generality to
$\G_m / F \,\, \times \,\, \G_m /F$ (see [Sch], p. 419), hence it 
does not extend in full generality to
algebraic tori $T$ with $dim \, T > 1.$ The phrase {\it full generality} 
in the last sentence means {\it for any} $P \in T(F)$ 
{\it and any subgroup} $\Lambda \subset T(F).$ Hence, as far as  full
generality for tori is concerned, the problem of detecting linear dependence by reduction maps has affirmative answer only for tori with $e=1.$

\medskip
\noindent\it The case of abelian varieties. 
\rm\medskip

\noindent
Let $E := E_{d}$ be the elliptic curve over $\Q$ given by the equation
$y^2 = x^3 -d^2 x.$ It has $CM$ by $\Z [i].$ It has been shown 
that the rank of $E_{d}(\Q)$ can reach $6$
see [RS], Table 2, p. 464. For example one can find in the 
Table loc. cit. that for $d = 34$ rank of $E_d (\Q)$ is $2,$ 
for $d = 1254$ rank of $E_d (\Q)$ is $3$ and for 
$d = 29274$ rank of $E_d (\Q)$ is $4$ (see [Wi]). 
Moreover for $d = 205015206$ the rank of $E_d (\Q)$ is $5$   
and for $d = 61471349610$ the rank of $E_d (\Q)$ is $6$ (see [Ro]).
>From now on we assume that the rank of $E_d (\Q)$ is at least 
$2.$ 

\noindent
Note that for every $d > 1$ the group $E_{p} (\F_p)$ does not have 
$p$ torsion for each $p \, {\not|}\,\, 2d.$ Indeed, for each $d > 1$
we have $E [2] \subset E(\Q).$ Hence by [Sil1], Prop. 3.1, p. 176
the group $E[2]$ injects into 
$E_{p} (\F_p)$ by the reduction map $r_p$ for every 
$p \, {\not|}\,\, 2d.$ Hence $4 \, | \, |E_{p} (\F_p)|$ for this $p.$  
On the other hand by the Theorem of Hasse we have $|E_{p} (F_p)| < p+1+2 \sqrt{p}$ which
implies that $|E_{p} (F_p)| < 4p$ for every $p \ge 3.$
This implies that $p$ does not divide $|E_{p} (\F_p)|$ for every $p \, {\not|}\,\, 2d.$ 

\noindent
Let us now consider the curve $E = E_d$ over $\Q (i).$ 
It is easy to observe that 
$\text{rank}_{\Z} \,\, E_{d} (\Q (i)) \, = \, 2 \text{rank}_{\Z} \,\, E_{d} (\Q ).$
Let $v$ denote a prime over $p$ for each $p \, {\not|}\,\, 2d.$ If $p$ splits completely
in $\Q (i) / \Q$ then $k_v = \F_p.$ In this case 
$E_v (k_v) = E_p (\F_p)$ and $E_v (k_v)$ does not have $p$ torsion.
If $p$ is inert in $\Q (i) / \Q,$ 
then by use of [Sil1] Theorem 4.1, c.f. p. 309 loc. cit.,
we observe that $E_v$ is supersingular, hence $E_v (k_v)$ does not have $p$-torsion by the theorem of Deuring [De] c.f. [Sil1], 
Theorem 3.1, p. 137.

\noindent
Note that $E (\C) \cong \C / \Z [i].$ Hence 
$E(\overline{\Q (i)})_{tor} \cong \Q(i) /\Z [i].$
On the other hand the reduction map gives a
natural isomorphism: 
$$E(\overline{\Q (i)})_{tor \not= p} \cong
 E_{v}(\overline{k_v})_{tor \not= p}.$$
Hence we can identify $E_v (k_v)$ with a subgroup
of $E[c] \cong {1 \over c} \Z [i] / \Z [i]$ for some 
$c \in \Z [i],$ \, $c \, {\not|} \, p.$
Note that in our case $E_{v}(k_v)$ is the fixed points of
the $Fr_v \in G(\overline{k_v} / k_v)$ acting on 
$E_{v}(\overline{k_v})_{tor \not= p}.$ 
Hence $E_{v}(k_v)$ is a cyclic $\Z [i]$-submodule of the cyclic
$\Z [i]$-module $E[c].$ So for each $p \, {\not|}\,\, 2d$ there is an element 
$\gamma (v) \in \Z [i]$
such that $E_{v}(k_v)$ is precisely the subgroup of $E[c]$ annihilated by 
multiplication by $\gamma (v).$ So for each $p \, {\not|}\,\, 2d$
we have $E_v (k_v) \cong {1 \over \gamma (v)} \Z [i] / \Z [i] \cong 
\Z [i] / \gamma (v).$ Hence $E_v (k_v)$ has a cyclic $\Z [i]$-module structure.

\noindent
We consider the  abelian surface $A_{d} := 
E_{d} \times E_{d} = E_{d}^2$ as defined over $\Q (i).$

\proclaim{Remark 6.1} 
{\rm For  abelian variety $A_{d}$ one has 
$e = 2 > dim_{\Q (i)} \, H_1 (E_{d} (\C); \, \Q) = 1.$ 
Hence $A$ is just beyond the range of abelian varieties 
considered in Theorem 5.1}
\endproclaim
\noindent In the proposition below we present a counterexample to the problem 
of detecting linear dependence for abelian varieties.

\proclaim{Proposition 6.2} There is a nontorsion point $P \in A_{d} (\Q (i))$ and 
a free $\Z [i]$-module $\Lambda \subset A_{d} (\Q (i))$ such that 
$P \notin \Lambda$ and
$r_{v}(P) \in r_{v} (\Lambda)$ for all primes $v\, \not | \, 2 d$ in $\Z [i].$
\endproclaim

\demo{Proof:}
By our assumption that rank of $E_d (\Q)$ is at least $2,$  we can find two points $Q_1, \, Q_2\,  \in \, E_d (\Q (i))$ such that they are independent over $\Z [i].$ 
Let $P, \, P_1, \, P_2, \,  P_3 \, \in \, A(\Q (i))$ be defined as follows:

$$P := \left[ \matrix
0\\
Q_1
\endmatrix \right], \,\, P_1 := \left[\matrix
Q_1\\
0
\endmatrix\right], \,\, P_2 := \left[ \matrix
Q_2\\
Q_1
\endmatrix \right], \,\, P_3 := \left[ \matrix
0\\
Q_2
\endmatrix \right].$$

\noindent Let $\Lambda := \Z [i] P_1 + \Z [i] P_2 + \Z [i] P_3 \subset A (\Q (i)).$
We observe that $\Lambda$ is free over $\Z [i]$ hence also free over $\Z.$ 
However $\Lambda$ is not free over $End_{\Q (i)} \, A = M_{2} \, (\Z [i]).$
Moreover it is clear that $P \notin \Lambda.$
\medskip  

\noindent
Let $\overline{Q_i} := r_v (Q_i)$ for $i = 1, 2,$ 
$\overline{P_i} := r_v (P_i)$ for $i = 1, 2, 3$ and 
$\overline{P} := r_v (P).$  We will prove that 
$r_v (P) \in r_v (\Lambda)$ for all $v$ of $\Z [i]$ over a prime 
$p \, {\not|}\,\, 2d.$ The equation 
$$\overline{P} = r_1 \overline{P_1} + r_2 \overline{P_2} + r_3 \overline{P_3}.$$
in $E_v (k_v) \times E_v (k_v)$ with $r_1, r_2, r_3 \in \Z [i]$
is equivalent to a system of equations in $E_v (k_v):$
$$r_1 \overline{Q_1} +  r_2 \overline{Q_2} = 0$$
$$r_2 \overline{Q_1} +  r_3 \overline{Q_2} =  \overline{Q_1}$$
Because $E_v (k_v) \cong \Z [i] / \gamma (v),$ there are elements 
$c_1, c_2 \in \Z [i]$ such that via this isomorphism we can make the following identifications $\overline{Q_1}\, = \, c_1 \mod \gamma (v)$ and 
$\overline{Q_2} \, = \, c_2 \mod \gamma (v).$ Hence the above system of equations is equivalent to the system of congruences in $\Z [i] / \gamma (v):$
$$r_1 c_1 +  r_2 c_2 \equiv 0 \mod \gamma (v) $$
$$r_2 c_1 +  r_3 c_2 \equiv   c_1 \mod \gamma (v).$$
If $c_1 \equiv 0 \mod \gamma (v) $ or $c_2 \equiv 0 \mod \gamma (v) $
then the last system of congruences trivially has a solution.
Hence assume that $c_1 \not\equiv 0 \mod \gamma (v) $ and 
$c_2 \not\equiv 0 \mod \gamma (v).$
Let $D := gcd (c_1, c_2).$ Then it is easy to check that
$$ gcd \, (c_{1}^{2} / D, \, c_2) \, = \, D$$
and since $D\, |\, c_1 $ 
it implies that the equation $r \,c_{1}^{2} / D \, + r_3 c_2 = c_1$
has a solution in $r, \, r_3 \, \in \, \Z [i].$ Putting 
$$r_1 := {-r c_2 \, \over \, D}, \quad r_2 := {r c_1 \, \over \, D}$$
we find out that numbers $r_1, \, r_2, \, r_3 \, \in \, \Z [i]$ satisfying
the above system of congruences. 
\qed\enddemo

\bigskip

\subhead  7. Detecting linear dependence via finite number of
reductions
\endsubhead
\medskip

Let $A/F$ be an abelian variety  defined over a
number field $F.$ Let
$$\beta_{H} \, :\, A(F) \otimes_{\Z} \R \, \times
A(F) \otimes_{\Z} \R \, \rightarrow\, \R$$
be the height pairing defined by the canonical height function
on $A$ [HS], [Sil2]. It is known loc. cit
that $\beta_{H}$ is positive definite, symmetric bilinear form.
Moreover if $R \in A(F)$ then  $\beta_{H} (R, R) = 0$ iff
$R$ is a torsion point.
\medskip

Let $P \in A (F)$ and let $\Lambda$ be a  subgroup of $A (F).$
Recall that $\Omega := c\, A (F).$
For our purposes, as explained in section 2, we will assume
that $\Lambda \subset \Omega.$ Let $r$ denote the rank of $\Lambda.$ 
Let
$P_1, \dots, P_r,\dots, P_s$ be such a $\Z$-basis of $\Omega$
that:
$$\Lambda = \Z d_1 P_1 + \dots +  \Z d_r P_r + \dots + \Z d_sP_s .
\tag{7.1}$$
where $d_i \in \Z \, \backslash \, \{0\}$ for $1 \leq i \leq r$ and
$d_i = 0$ for $i > r.$ For any $P \in A(F)$ we can write
$$cP = \sum_{i=1}^s  n_i P_i \tag{7.2}$$
and we get
$$c^2 \, \beta_{H}(P, P) = \sum_{i, j} n_i n_j
\beta_{H}(P_{i}, P_{j}). \tag{7.3}$$
Since $\beta_{H}(P, P) > 0$ and $\beta_{H}$ is positive definite,
there is a constant $C$ which depends only on the points $P,
P_{1}, \dots, P_{s}$ such that

$$|n_i| \leq C, \,\, \text{for all}\,\, 1 \leq i \leq s. \tag{7.4}$$
Hence if $P \in \Lambda$ then
$P = \sum_{i = 1}^r k_i d_i P_{i}$
for some $k_1, \dots, k_r \in \Z.$ It follows that $|d_i k_i| \leq C,$
so $|k_i| \leq {C \over d_i} \leq C$ for each $1 \leq i \leq r.$
Hence there is only a finite number, not bigger than $(2 C + 1)^r,$ of
tuples $(n_1, \dots, n_r)$ to check to determine if $P \in \Lambda.$

We will apply the estimation of coefficients (7.4)
obtained by application of the height pairing in the proof of
Theorem 7.7.

\proclaim{Theorem  7.5} Let $A = A_1 \times \dots\times A_t$ 
be a product of simple, pairwise nonisogenous abelian varieties.
Let $l$ be a prime number and let
$Q_{ij} \in A_i (L)$ for $1 \leq j \leq
r_i$ be independent over ${\cal R}_i$ for each $1 \leq i \leq t.$
Let $L/F$ be a finite extension and $L_{l^m}:=L(A[l^m]).$ Let $k$ be
a natural number such that the image of
${\overline{\rho}}_{{l}^{k+1}}: G_{L_{l^k}}\rightarrow GL_{{\Bbb
Z}/{l^{k+1}}}(A[l^{k+1}])$ contains a nontrivial homothety.  Let $d$
be a discriminant of $L_{l^{k+1}}({\frac{1}{l^k}}{\Gamma})/{\Q}.$
There are effectively computable constants $b_{1}$ and $b_{2}$ such
that $r_w (Q_{ij}) = 0$ in $A_{i\, w} (k_w)_l$ for all $1 \leq j
\leq r_i$ and $1 \leq i \leq t$ for some prime $w$ of ${\cal O}_L$
such that $N_{L/{\Q}}(w)\le b_{1}d^{b_{2}}.$
\endproclaim

\demo{Proof} We argue in the same way as in the proof of
Theorem 3.3 but instead of using classical Chebotarev's theorem
we use the effective Chebotarev's theorem [LO] p. 416.
\qed\enddemo

\proclaim{Theorem  7.6} Let $A = A_1 \times \dots\times A_t$ 
be a product of simple, pairwise nonisogenous abelian varieties.
Let $l$ be a prime number. Let $m \in \N
\cup \{0\}$ for all $1 \leq j \leq r_i$ and $1 \leq i \leq t.$ Let
$L/F$ be a finite extension and let $P_{ij} \in A_i (L)$ be
independent over ${\cal R}_i$ and let $T_{ij} \in A_i [l^m]$ be
aribitrary torsion elements for all $1 \leq j \leq r_i$ and $1 \leq
i \leq t.$  Let $k\ge m$ be a natural number such that the image of
${\overline{\rho}}_{{l}^{k+1}}: G_{L_{l^k}}\rightarrow GL_{{\Bbb
Z}/{l^{k+1}}}(A[l^{k+1}])$ contains a nontrivial homothety.  Let $d$
be a discriminant of $L_{l^{k+1}}({\frac{1}{l^k}}{\Gamma})/{\Q}.$
There are effectively computable constants $b_{1}$ and $b_{2}$
and there is a prime $w$ of ${\cal O}_L$ such
that $N_{L/{\Q}}(w)\leq b_{1}d^{b_{2}}$ and
$$r_{w^{\prime}} (T_{ij}) = r_w (P_{ij}) \,\,\,
\text{in}\,\,\, A_{i, w} (k_w)_l$$ for all $1 \leq j \leq r_i$ and
$1 \leq i \leq t,$ where $w^{\prime}$ is a prime in ${\cal O}_{L(A_i
[l^m])}$ over $w$ and $r_{w^{\prime}}\, :\, A_{i} (L(A_i [l^m]))
\rightarrow A_{i, w} (k_{w^{\prime}})$ is the reduction map.
\endproclaim
\demo{Proof} Follows immediately from Theorem 7.5 in the same way as
the Theorem 3.6 follows from Theorem 3.3.
\qed\enddemo
\medskip

\proclaim{Theorem 7.7} Let $A/F$ satisfy the hypotheses of Theorem
5.1. Let $P \in A(F)$ and let $\Lambda$ be a subgroup of $A(F).$
There is a finite set $S^{fin}$ of primes $v$ of ${\cal O}_F,$
depending on $A, P, \Lambda $ and the basis $P_1, \dots, P_s$ such
that the following condition  holds: if $r_v (P) \in r_v (\Lambda)$
for all $v \in S^{fin}$ then $P \in \Lambda + A(F)_{tor}.$

\noindent
Hence if $A(F)_{tor} \subset \Lambda$
then the following conditions are equivalent:
\roster
\item[1]\,\, $P \in \Lambda$
\item[2]\,\, $r_v (P) \in r_v (\Lambda)$ for all
$v \in S^{fin}.$
\endroster
\endproclaim

\demo{Proof} To construct the set $S^{fin}$ we will carefully  analyze the
proof of Theorem 5.1.  The finitness of $S^{fin}$ will follow by
application  of both the  canonical height function and the Theorem
of Lagarias and Odlyzko [LO] p. 416.
By explanation similar to that in section 2 we can assume, 
that $P\in\Omega$ and $\Lambda\subset\Omega$.
Consider the projections ${\pi}_{i}:{\Omega}\rightarrow {\Bbb Z}, $
${\pi}_j (R)= \mu_{j}, \, j=1,\dots ,s$ for $R={\sum}_{j=1}^{n}{\mu}_{j}P_{j}.$
In the same way as in the proof of the Theorem 5.1 construct for each
${\pi_j}$ the homomorphism
${\widetilde\pi}_{j} \in  Hom_{\M_{\text{e}} (\DD)}
(\Omega \otimes_{\Z} \Q, \,\, \M_{\text{e}} (\DD))$
such that $tr( \widetilde{\pi_i}) = \pi_i .$
Simiarly as in the proof of Theorem 5.1 we construct
 the maps:
${\widetilde s}_j,$
${\widetilde \pi (i)}_j,$ ${\widetilde s (i)}_j,$ where
${\widetilde \pi}_j = \prod_{i=1}^t {\widetilde \pi (i)}_j,$ \,
${\widetilde s}_j = \prod_{i=1}^t {\widetilde s (i)}_j.$ Moreover
$ \, \text{Ker} \, \widetilde{\pi} =
\prod_{i}^t \,\, \text{Ker} \,\, \widetilde{\pi (i)}.$
Then we construct the number $M_{0, j}$ and
the lattice
$$\widetilde{\Gamma (i)_j} := \sum_{k=1}^{k_{i, j}}\,\,
\,\, M_{0, j} \,\,  {\cal R}_i \widetilde{\omega_k (i)_{j}}
+ \sum_{k=k_{i, j}+1}^{u_{i, j}} {\cal R}_i \,\,
\widetilde{\omega_k (i)_j} \subset {\Omega}_{i} \otimes_{\Z} {\cal O}_K$$
and then the lattice $\widetilde{\Gamma_j} :=
\bigoplus_{i=1}^t \,\widetilde{\Gamma (i)_j}.$
Then we define numbers $M_{2, j}$ and $M_{3, j}$ such that
$M_{2, j} := [\Omega \otimes_{\Z} {\cal O}_{K}\, :\, \widetilde{\Gamma_j}]$
and  $M_{3, j} :=
[\widetilde{\Gamma_j}\, :\, M_{2, j}\, \Omega \otimes_{\Z} {\cal O}_{K}].$
For $n_j \not= 0$ in decomposition of $P$ in formula (2.2)
we consider every $l | n_j$ and every $\lambda | l$ and consider
the ramification index $\epsilon_{j, \lambda}$
of $\lambda$ over $l.$ Next we define $m_{1,j, \lambda}$ such that
$\lambda^{m_{1,j, \lambda}} \, ||\, n_{j}.$
We put $m_{2,j, \lambda} := m_{1,j, \lambda} + 1$ and
$m_{j, \lambda} := m_{2,j, \lambda} + 1.$
Following the proof of Theorem 5.1 we also construct
the constant $M_{4}$ which is clearly independent of $j.$
We define the nonnegative integers $m_{0,j}, m_{3,j}, m_4$
with the property $l^{m_{0,j}}\, || \, M_{0, j},$ \,
$l^{m_{3,j}}\, || \, M_{3, j}$
and $l^{m_{4}} \, || \, M_{4}.$ Put 
$m_{j,l} := \, \max_{\lambda\, |\, l}\,\, m_{j, \lambda},$ and 
$\epsilon_{j,l} := \, \max_{\lambda\, |\, l}\,\, \epsilon_{j, \lambda}.$
Now, we choose the
number $n_{j, l}$ in  such a way  that
the image of the representation
$${\overline{\rho}}_{{l}^{n_{j, l} + 1}}: 
G_{L_{l^{n_{j, l}}}}\rightarrow GL_{{\Bbb
Z}/{l^{n_{j, l} + 1}}}(A[l^{n_{j, l} + 1}])$$ 
contains a nontrivial homothety 
and $n_{j, l} >
\epsilon_{j, l}\, m_{0,j} + \epsilon_{j, l} \,
m_{4} + m_{j, l} + \epsilon_{j, l} \, m_{3,j}.$ The last inequality 
guaranties that 
$ \epsilon_{j, \lambda} n_{j, l}  >
\epsilon_{j, \lambda}\, m_{0,j} + \epsilon_{j, \lambda} \,
m_{4} + m_{j, \lambda} + \epsilon_{j, \lambda} \, m_{3,j}.$
Eventually, we construct for each
$1 \leq j \leq s$ and for each prime number
$l \, |\, \pi_j (P)$ the number field
$L_{j, l} := F(r(A)[l^{n_{j} + 1}], {1 \over l^{n_j}} \, 
\widetilde{\Gamma_j}),$
where $r (A)$ is the radical of $A$ defined in section 2.
Observe that there are only finite number of primes $l$
considered above by the estimation of coefficients (7.4).
By the Theorem of Lagarias and Odlyzko [LO] p. 416 there are
effectively computable constants $b_1$ and $b_2$ such that
every element $\sigma \in G(L_{j, l} / F)$
is equal to a Frobenius element $Fr_v \in G(L_{j, l} / F)$
for a prime $v$ of ${\cal O}_F$ such that
$N_{F/\Q} (v) \leq b_1 d_{L_{j, l}}^{b_2}.$
Now for every $j$ such that $n_j = \pi_j (P) \not= 0$ let
$$S_{j, l}^{fin} := \{v\, :\,
N_{F/\Q} (v) \leq b_1 d_{L_{j, l}}^{b_2} \,\, \text{and}\,\,
v\,\, \text{is of good reduction for}\,\, A \},$$
$$S_{j}^{fin} := \bigcup_{l | n_j}\,\, S_{j, l}^{fin}.$$
Then we define
$$S^{fin} \,\,:= \,\,  \bigcup_{1 \leq j \leq s, n_j \not= 0}
\, S_{j}^{fin}.$$
It is enough to prove that for the set $S^{fin}$ condition (2)
implies (1). Indeed, if (1) does not hold then in the same way as in
he proof of the Theorem 5.1 there is $1 \leq j_0 \leq s$ such that
$P \notin \Lambda \otimes_{\Z} {\cal O}_{\lambda}$ for some $l$ and
$\lambda\, | \, l$ such that $\lambda^{m_{1, j_0, \lambda}}\, ||\, n_{j_0}$
and $\lambda^{m_{2, j_0, \lambda}}\, |\, d_{j_0}$ for natural numbers $m_{1,
j_0, \lambda} < m_{2, j_0, \lambda} = m_{1,
j_0, \lambda} + 1.$ As in the proof of Theorem 5.1 this leads to the 
investigation of a homomorphism $\pi_{j_0}$ of $\Z$-modules and now the 
proof follows the lines of the proof of Theorem 5.1. 
Of course, the choice of prime $w$
in ${\cal O}_{F(r(A)[l^{n_{j_0}}])}$ is done now by virtue of 
Theorem 7.6. So it is
clear by the definition of $S_{j_0}^{fin}$ that such a prime $w$ can
be chosen over a prime $v \in S_{j_0}^{fin}.$ Hence in the same way
as in the proof of Theorem 5.1 we are led to a contradiction.
\qed\enddemo

\proclaim{Remark 7.8} \rm{The problem with an effective algorithm
for finding $S^{fin}$ comes from the lack of an effective algorithm
for finding the $\Z$-basis of $A(F)/A(F)_{tor}.$ See [HS] p. 457-465
for the explanation of the obstructions for an effective algorithm for
finding the $\Z$-basis of $A(F)/A(F)_{tor}.$}
\endproclaim

\proclaim{Remark 7.9}
\rm{For a given abelian variety $A/F$, in general, there is
no finite set $S^{fin}$ of primes of good reduction, that depends only
on $A,$ such that for any
$P \in A(F)$ and any subgroup $\Lambda \in A(F)$ the condition
$r_v (P) \in r_v (\Lambda)$ for
all $v \in S^{fin}$ implies $P \in \Lambda + A(F)_{tor}.$
Indeed, take any simple abelian variety $A$ with
$End_{\overline{F}} (A) = \Z$
and rank of $A(F)$ over $\Z$ at least $2.$
Take two nontorsion points $P^{\prime}, Q^{\prime} \in A(F),$
linearly independent
over $\Z.$ For any natural number $M$ consider
the finite set $S_{M}$ of primes $v$ of ${\cal O}_F$ of good reduction
for $A/F$ which are over rational primes $p \leq M.$
Take a natural number $n$ divisible by
$\prod_{v \in S_M} \,\, |A_{v} (k_v)|.$ Taking $P := n P^{\prime}$
and $\Lambda := n \Z\, Q^{\prime}$ we observe that
$r_v (P) = 0 = r_{v} (\Lambda)$ for all $v \in S_{M}$
but by construction $P \notin \Lambda + A(F)_{tor}.$}
\endproclaim
\bigskip

\noindent
{\it Acknowledgements}:\quad
The authors would like to thank A. Schinzel for pointing out
the example in his paper [Sch] on page 419. The authors would like to
thank K. Rubin for pointing out important properties of CM elliptic curves  
used in section 6 of this paper and for a number of important suggestions. 
The research was partially financed by
a research grant of the Polish Ministry of Science
and Education.

\enddocument